\documentclass[3p,times]{elsarticle}

\usepackage{amsmath}

\usepackage{amssymb}
\usepackage[figuresright]{rotating}

\usepackage[english]{babel}
\begin{document}
	
 \begin{frontmatter}

  \title{Umbral Methods and Harmonic Numbers}

	\author[Enea]{G. Dattoli \corref{cor}}
	\ead{giuseppe.dattoli@enea.it}
		
	\author[LaS]{B. Germano}
	\ead{bruna.germano@sbai.uniroma1.it}
		
	\author[Enea,Unict]{S. Licciardi }
	\ead{silvia.licciardi@dmi.unict.it}
		
	\author[LaS]{M. R. Martinelli}
	\ead{martinelli@dmmm.uniroma1.it}

	 \address[Enea]{ENEA - Frascati Research Center, Via Enrico Fermi 45, 00044, Frascati, Rome, Italy}
	 \cortext[cor]{Corresponding author}
	 \address[LaS]{University of Rome, La Sapienza, Department of Methods and Mathematic Models for Applied Sciences, Via A. Scarpa, 14, 00161 Rome, Italy}
	 \address[Unict]{University of Catania, Department of Mathematics, Viale Doria, 6, 95125 Catania, Italy}
	 \address[LaS]{University of Rome, La Sapienza, Department of Methods and Mathematic Models for Applied Sciences, Via A. Scarpa, 14, 00161 Rome, Italy}
		
\begin{abstract}
The theory of harmonic based function is discussed here within the framework of umbral operational methods. We derive a number of results based on elementary notions relying on the properties of Gaussian integrals.
\end{abstract}
		
		
 \end{frontmatter}
	
\section{Introduction}

 Methods employing the concepts and the formalism of umbral calculus have been exploited in \cite{Dattoli} to conjecture the existence of generating functions involving Harmonic Numbers \cite{Sondow}. The conjectures put forward in \cite{Dattoli} have been proven in \cite{Coffee}-\cite{Cvijovi}, further elaborated in subsequent papers \cite{Mezo} and generalized to Hyper-Harmonic Numbers in \cite{Conway}.

\noindent In this note we use the same point of view of \cite{Dattoli} , by discussing the possibility of exploiting the formalism developed therein in a wider context.\\

\noindent We accordingly remind that harmonic numbers are defined as
\begin{equation} \label{GrindEQ__1_} 
h_{n} =\sum _{r=1}^{n}\frac{1}{r}   
\end{equation} 
It is furthermore evident that the integral representation for this family of numbers can be derived using a standard procedure, we first note that

\begin{equation}\label{int}
h_{n} =\sum _{r=1}^{n}\int _{0}^{\infty }e^{-s\, r}   ds
\end{equation}
thereby getting \cite{Rochowicz}

\begin{equation}\label{gau}
h_{n} =\int _{-\infty }^{0}\frac{e^{(n+1)\, \xi } -e^{\xi } }{e^{\xi } -1}  d\xi =\int _{0}^{1}\frac{1-x^{n} }{1-x} dx 
\end{equation}
 after interchanging summation and integral signs. 

\noindent The definition in eq. \eqref{gau} can be extended to non-integer values of $n$ and, therefore, it can be exploited as an alternative definition holding for $n$  (not-necessarily) a positive real. \\

 We define now the umbral operator 
\begin{equation}\begin{split} \label{GrindEQ__3_} 
& \hat{h}^{n} =h_{n} , \\ 
& h_{0} =1 
\end{split}\end{equation} 
with the property

\begin{equation} \label{GrindEQ__4_} 
\hat{h}^{n} \hat{h}^{m} =\hat{h}^{n+m}  
\end{equation} 
and introduce the Harmonic Based Exponential Function \textit{(HBEF)}
\begin{equation} \label{GrindEQ__5_} 
{}_{h} e(x)=e^{\hat{h}\, x} =1+\sum _{n=1}^{\infty }\frac{h_{n} }{n!} \,  x^{n}  
\end{equation} 
which, as already discussed in \cite{Dattoli}, has quite remarkable properties.\\

\noindent The relevant derivatives can accordingly  be expressed as (see the concluding part of the paper for further comments)

\begin{equation} \begin{split}\label{GrindEQ__6_} 
& \left(\frac{d}{dx} \right)^{m} {}_{h} e(x)={}_{h} e(x,m)=\hat{h}^{m} e^{\hat{h}\, x} =h_{m} +\sum _{n=1}^{\infty }\frac{h_{n+m} }{n!} \,  x^{n} , \\ 
& \left(\frac{d}{dx} \right)^{m} {}_{h} e(x,k)={}_{h} e(x,k+m) 
\end{split}\end{equation} 
We use the previous definition to derive the following integral 

\begin{equation} \label{GrindEQ__7_} 
\int _{0}^{\infty }{}_{h} e (-\alpha \, x)\, e^{-x} dx=\int _{0}^{\infty }e^{-(\alpha \, \hat{h}+1)\, x}  dx=\frac{1}{\alpha \, \hat{h}+1}  
\end{equation} 
It is evident that by expanding the umbral function on the r h s of eq. \eqref{GrindEQ__7_}, we obtain 
\begin{equation} \label{GrindEQ__8_} 
\frac{1}{\alpha \, \hat{h}+1} =1+\sum _{s=1}^{\infty }(-1)^{s}  \alpha ^{s} h_{s}  
\end{equation} 
an expected conclusion, achievable by direct integration, underscored  here to stress the consistency of the procedure.\\

\noindent A further interesting example comes from the following ``Gaussian'' integral

\begin{equation}\label{krak}
\int _{-\infty }^{\infty }{}_{h} e (-\alpha \, x)\, e^{-x^{2} } dx=\int _{-\infty}^{\infty }e^{-(\alpha \, \hat{h}\, x+x^{2} )\, }  dx=\sqrt{\pi } e^{\frac{\alpha ^{2} \, \hat{h}^{2} }{4} }
\end{equation} 
 The last term in eq. \eqref{krak} has been obtained by treating $\hat{h}$ as an ordinary algebraic quantity and then by applying the standard rules of the Gaussian integration


\begin{equation}
e^{\frac{\hat{h}^2 \alpha^2}{4}}=
{}_{h^2}e\left( \dfrac{\alpha^2}{4}\right)= 1+\sum _{r=1}^{\infty }\frac{h_{2r} }{r!}  \left(\frac{\alpha }{2} \right)^{2\, r} 
\end{equation}

\noindent Let us now consider the following slightly more elaborated example, involving the integration of two ``Gaussians'', namely the ordinary case and its HBEF analogous

\begin{equation} \label{GrindEQ__10_} 
\int _{-\infty }^{\infty }{}_{h} e (-\, \alpha \, x^{2} )\, e^{-\, x^{2} } dx=\int _{-\infty }^{\infty }e^{-(\, \hat{h}\, \alpha +\, 1)x^{2} \, }  dx=\sqrt{\frac{\pi }{1+\alpha \, \hat{h}} }  
\end{equation} 
This last result, obtained after applying elementary rules, can be worded as it follows: the integral in eq. \eqref{GrindEQ__10_} depends on the operator function on its r.h.s., for which we should provide a computational meaning. The use of the Newton binomial yields

\begin{equation}\begin{split} \label{GrindEQ__11_} 
& \sqrt{\frac{\pi }{1+\alpha \hat{h} } } =\sqrt{\pi } \sum _{r=0}^{\infty }\binom{-\frac{1}{2}}{r} \, \left(\alpha \, \hat{h}\right)^{\, r} =\sqrt{\pi } \left(1+\sqrt{\pi } \sum _{r=1}^{\infty }\frac{\alpha ^{r} h_{r} }{\Gamma \left(\frac{1}{2} -r\right)\, r!}  \, \right), \\ 
& \left|\alpha \right|<1
 \end{split} \end{equation} 
and the correctness of this conclusion has been confirmed by the numerical check.\\

\noindent It is evident that the examples we have provided show that the use of concepts borrowed from umbral theory offer a fairly powerful tool to deal with the ``harmonic based'' functions.\\

 The next step we will touch in this paper is to check whether non integer forms of harmonic numbers make any sense.\\ 

 We consider indeed the following function
 
\begin{equation} \label{GrindEQ__12_} 
{}_{\sqrt{h} } e(x)=e^{\hat{h}^{\frac{1}{2} } \, x} =1+\sum _{n=1}^{\infty }\frac{\left(\sqrt{\hat{h}} \, x\right)^{n} }{n!} \,  =1+\sum _{n=1}^{\infty }\frac{h_{n/2} }{n!} \,  \left(x\right)^{n}  
\end{equation} 
and the integral

\begin{equation}\begin{split} \label{GrindEQ__13_} 
& \int _{-\infty }^{+\infty }{}_{\sqrt{h} } e (\alpha \, x)\, e^{-x^{2} } dx=\int _{-\infty }^{+\infty }e^{\hat{h}^{\frac{1}{2} } \, \alpha \, x-x^{2} }  dx= \\ 
& =\sqrt{\pi } e^{\hat{h}\, \left(\frac{\alpha }{2} \right)^{2} \, } =\sqrt{\pi } {}_{h} e\left(\left(\frac{\alpha }{2} \right)^{2} \right)  
\end{split}\end{equation} 
To this aim we remind the following identity from Laplace transform theory \cite{Doetsch}

\begin{equation}\begin{split} \label{GrindEQ__14_} 
& e^{-p^{\frac{1}{2} } \, x} =\int _{0}^{\infty }e^{-p\, \eta \, x^{2} }  \, g_{\frac{1}{2} } (\eta )\, d\eta  \\ 
& g_{\frac{1}{2} } (\eta )\, =\frac{1}{2\, \sqrt{\pi \eta ^{3} } } e^{-\frac{1}{4\, \eta } } \end{split}\end{equation} 
The use of eq. \eqref{GrindEQ__12_} allows to write the identity

\begin{equation}\begin{split} \label{GrindEQ__15_} 
& {}_{\sqrt{h} } e(-x)=\int _{0}^{\infty }{}_{h} e (-\eta \, x^{2} )\, g_{\frac{1}{2} } (\eta )\, d\eta , \\
&  g_{\frac{1}{2} } (\eta )\, =\frac{1}{2\, \sqrt{\pi \eta ^{3} } } e^{-\frac{1}{4\, \eta } } 
\end{split}\end{equation} 
the numerical check has, in both cases, confirmed the correctness of the ansatz.\\
 The possibility of defining ${}_{\sqrt[{k}]{h} } e(x)$ will be discussed elsewhere.\\

 We go back to eq. \eqref{GrindEQ__6_} and write the first derivative of the HBEF as
\begin{equation} \label{GrindEQ__16_} 
_{h} e(x,1)=1+\sum _{n=1}^{\infty }\frac{h_{n+1} }{n!} \,  x^{n}  
\end{equation} 
By taking into account that $h_{n+1} =h_{n} +\frac{1}{n+1} $ we end up with the following differential equation defining the function $_{h} e(x)$

\begin{equation} \label{GrindEQ__17_} 
y'=y+\frac{e^{x} -1-x}{x}  
\end{equation} 
The relevant solution reads

\begin{equation}\begin{split} \label{GrindEQ__18_} 
& {}_{h} e(x)=1+e^{z} \left(\ln (x)+E_{1} (x)+\gamma \right), \\ 
& E_{1} (x)=\int _{x}^{\infty }\frac{e^{-t} }{t}  dt, \\ 
& \left(\ln (x)+E_{1} (x)+\gamma \right)=-\sum _{n=1}^{\infty }\frac{(-x)^{n} }{n\, n!}  , \\ 
& \gamma \equiv  Euler- Mascheroni- constant  
\end{split}\end{equation} 
which is the generating function of harmonic numbers originally derived by Gosper (see \cite{Sondow}).\\

\noindent By iterating the previous procedure we find the following general recurrence

\begin{equation}\label{rec}
_{h} e(x,m)={}_{h} e(x)+\sum _{r=0}^{m-1}\left(\frac{d}{dx} \right)^{r}  \frac{e^{x} -1-x}{x}
\end{equation}     

The harmonic polynomials

\begin{equation} \label{GrindEQ__20_} 
h_{n} (x)=(\hat{h}+x)^{n} =\sum _{s=0}^{n}\binom{n}{s}\, h_{s} x^{n-s}   
\end{equation} 
are easily shown to be linked to the HBEF by means of the generating function

\begin{equation} \label{GrindEQ__21_} 
\sum _{n=0}^{\infty }\frac{t^{n} }{n!}  h_{n} (x)=e^{x\, t} {}_{h} e(t) 
\end{equation} 
They belong to the family of App\'el polynomials and satisfy the recurrences

\begin{equation}\begin{split} \label{GrindEQ__22_} 
& \frac{d}{dx} h_{n} (x)=n\, h_{n-1} (x), \\ 
& h_{n+1} (x)=(x+1)\, h_{n} (x)+f_{n} (x), \\ 
& f_{n} (x)=\sum _{s=0}^{n}\frac{n!}{s!\, (n-s)!}  \frac{x^{n-s} }{s+1} =\int _{0}^{1}(x+y)^{n}  dy
\end{split} \end{equation} 
along with the identity, proved by iteration

\begin{equation}
h_{n} (-1)=(-1)^n \left( 1-\frac{1}{n}\right) 
\end{equation}
and
\begin{equation} h_{n} =\sum _{s=0}^{n}\binom{n}{s}\, h_{s}  (-1)
 \end{equation}                
the harmonic Hermite polynomials (touched on in ref. \cite{Dattoli}-\cite{Coffee}-\cite{Zhukovsky}) can also be written as
 
\begin{equation}\begin{split} \label{GrindEQ__24_} 
& \sum _{n=0}^{\infty }\frac{t^{n} }{n!}\;  {}_{h} H_{n} (x)=e^{x\, t} {}_{h} e(t^{2} ), \\ 
& H_n (x, \hat{h})= e^{\hat{h}\partial_x ^2}x^n =n!\, \sum _{r=0}^{\lfloor\frac{n}{2} \rfloor}\frac{x^{n-2\; r} \hat{h}^{r} }{(n-2\, r)!\, r!},\\
& {}_{h} H_{n} (x)=n!\, \sum _{r=0}^{\lfloor\frac{n}{2} \rfloor}\frac{x^{n-2\; r} h_{r} }{(n-2\, r)!\, r!}  
 \end{split} \end{equation} 
The use of the methods put forward in \cite{GDattoli} yields the recurrences

\begin{equation}\begin{split} \label{GrindEQ__25_} 
& \frac{d}{dx} {}_{h} H_{n} (x)=n\, {}_{h} H_{n-1} (x) \\ 
& {}_{h} H_{n+1} (x)=\left(x+2\, \hat{h}\, \frac{d}{dx} \right)\, {}_{h} H_{n} (x)=\left(x+2\frac{d}{dx} \right)\, {}_{h} H_{n} (x)+2\, \alpha '_{n} (x) \\ 
& \alpha _{n} (x)=n!\sum _{s=1}^{\lfloor\frac{n}{2} \rfloor}\frac{x^{n-2\, s} }{s!\, (n-2s)!}  \frac{1}{(s+1)} =\int _{0}^{1}(H_{n}  (x,\, y)-x^n)dy,\,  \\ 
& \alpha '_{n} (x)=\frac{d}{dx} \alpha _{n} (x)
 \end{split} \end{equation} 
Thus getting for the corresponding differential equation the following non homogeneous $ODE$

\begin{equation} \label{GrindEQ__26_} 
\left(x\frac{d}{dx} +2\left(\frac{d}{dx} \right)^{2} \right)\, {}_{h} H_{n} (x)=n \;{}_{h}H_{n} (x)-2\, \alpha ''_{n} (x) 
\end{equation}
 
Before closing the paper, we want to stress the possibility of extending the present procedure to the truncated exponential numbers, namely
\begin{equation} \label{GrindEQ__27_} 
e_{n} =\sum _{r=0}^{n}\frac{1}{r!}   
\end{equation} 
The relevant integral representation writes \cite{Marinelli}

\begin{equation} \label{GrindEQ__28_} 
e_{\alpha } =\frac{1}{\Gamma (\alpha +1)} \int _{0}^{\infty }e^{-s}  (1+s)^{\alpha } ds 
\end{equation} 
which holds for non-integer real values of $\alpha $ too. For example we find

\begin{equation}\label{Gammamez}
e_{-\frac{1}{2} } =\frac{e}{\sqrt{\pi } } \Gamma \left( \frac{1}{2},1 \right) 
\end{equation}               
with $\Gamma \left( 1,\frac{1}{2} \right) $ being the truncated Gamma function. According to the previous  discussion and to eq. \eqref{Gammamez}, setting $\hat{e}^\alpha \leftrightarrow e_{\alpha}$, we also find that
\begin{equation}\begin{split} \label{GrindEQ__30_} 
& \int _{-\infty }^{+\infty }e^{-\hat{e}\, x^{2} }  dx=\sqrt{\pi } e_{-\frac{1}{2} } , \\ 
& e^{-\hat{e}\, x^{2} } =\sum _{r=0}^{\infty }(-1)^{r}  \frac{e_{r} }{r!} x^{2\, r} 
 \end{split} \end{equation} 
This last identity is a further proof that the implications offered by the topics treated in this paper are fairly interesting and deserve further and more detailed investigation, which will be more accurately treated elsewhere.\\

\end{document}